\documentclass{amsart}
\usepackage{amsmath,amsfonts,amssymb}
\usepackage{graphicx}
\usepackage{color}

\def\Z#1{{\mathbb Z}^{#1}}
\def\Q#1{{\mathbb Q}^{#1}}
\def\R#1{{\mathbb R}^{#1}}
\def\C#1{{\mathbb C}^{#1}}
\def\T#1{{\mathbb T}^{#1}}

\def\noz{m}
\def\noi{n}

\def\cH{{\mathcal H}}
\def\cL{{\mathcal L}}
\def\cM{{\mathcal M}}

\newcommand{\SL}{\operatorname{SL}(2,{\R{}})}
\newcommand{\GL}{\operatorname{GL}(2,{\R{}})}

\newtheorem*{Theorem}{Theorem}
\newtheorem*{MagicWand}{Magic Wand Theorem}

\begin{document}

\begin{picture}(0,0)(0,0)
\put(0,60){English translation of}
\put(0,50){\textit{``Le  th\'eor\`eme  de  la  baguette  magique  de  A.~Eskin  et
M.~Mirzakhani''},}
\put(0,40){Gazette  des  math\'ematiciens,  \textbf{142} (2014) 39--54.}
\end{picture}

\title{The Magic Wand Theorem of A.~Eskin and M.~Mirzakhani}

\author{Anton Zorich}
\address{
Institut de math\'ematiques de Jussieu --- Paris Rive Gauche,
Institut Universitaire de France}
\email{Anton.Zorich@gmail.com}

\maketitle


\vspace*{-20pt}

\includegraphics{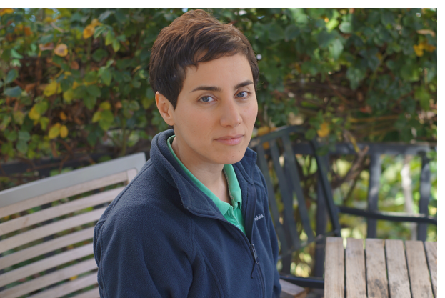}
\vspace*{185pt}
\begin{picture}(0,0)
\put(187,16){\tiny Courtesy of Maryam Mirzakhani}
\end{picture}

On  August  13,  2014  (the  openning  day  of  ICM  at Seoul) Maryam
Mirzakhani  received  the  Fields Medal \textit{``for her outstanding
contributions  to  the  dynamics and geometry of Riemann surfaces and
their  moduli  spaces''}  becoming  the first woman to win the Fields
Prize. Citing the ICM laudation~\cite{ICM},

\textit{``Maryam  Mirzakhani has made stunning advances in the theory
of  Riemann  surfaces and their moduli spaces, and led the way to new
frontiers  in  this  area.  Her insights have integrated methods from
diverse  fields, such as algebraic geometry, topology and probability
theory.}

\textit{In  hyperbolic  geometry,  Mirzakhani  established asymptotic
formulas  and statistics for the number of simple closed geodesics on
a Riemann surface of genus $g$. She next used these results to give a
new and completely unexpected proof of Witten's conjecture, a formula
for  characteristic classes for the moduli spaces of Riemann surfaces
with marked points.}

\textit{In  dynamics,  she  found  a remarkable new construction that
bridges  the  holomorphic and symplectic aspects of moduli space, and
used  it  to  show  that  Thurston's  earthquake  flow is ergodic and
mixing.}

\textit{Most  recently,  in  the  complex  realm,  Mirzakhani and her
coworkers produced the long sought-after proof of the conjecture that
--–  while  the  closure  of a real geodesic in moduli space can be a
fractal  cobweb,  defying classification --– the closure of a complex
geodesic is always an algebraic subvariety.}

\textit{Her work has revealed that the rigidity theory of homogeneous
spaces  (developed  by  Margulis,  Ratner  and others) has a definite
resonance  in the highly inhomogeneous, but equally fundamental realm
of moduli spaces, where many developments are still unfolding.''}

We   start   this   article   with   a  short  biographical  note  of
Maryam~Mirzakhani    (for   more   details   see   excellent   online
article~\cite{Quanta}  of  Erica~Klarreich).  Then  we try to give an
idea  of  her  achievements  focusing on one fundamental result as an
example.   To   put   this   result   into   context  we  present  in
section~\ref{s:phisics:2:foliations}   several  natural  problems  in
physics   which   lead   to   measured  foliations  on  surfaces.  In
section~\ref{s:Translation:surfaces}   we   describe   the  extremely
elementary  and  convenient  language of translation surfaces used to
work  with  measured  foliations. In section~\ref{s:modulis:space} we
provide background material describing how dynamics on one individual
translation  surface leads to dynamics in a multidimensional space of
translation  surfaces.  In  section~\ref{s:almost:all} we discuss why
ergodic theory is usually absolutely powerless in describing specific
trajectories. And, finally, in section~\ref{s:Theorem} we present the
Theorems   of   A.~Eskin   and   M.~Mirzakhani,   and   of  A.~Eskin,
M.~Mirzakhani, A.~Mohammadi, which serve as a Magic Wand for numerous
applications (the choice of applications in this exposition is mostly
based   on   simplicity   of   presentation   rather  than  on  their
significance).

This article is closely related to the part of the paper of P.~Hubert
and  R.~Krikorian  in  the  current  issue  dedicated  to  results of
Artur~Avila      on      dynamics     in     Teichm\"uller     space,
see~\cite{Hubert:Krikorian}.  The reader might find it interesting to
read  both  companion  papers:  they  are intended to complement each
other.

For    more    ample    presentation   of   mathematical   works   of
Maryam~Mirzakhani   see   the  short  paper~\cite{McMullen:laudition}
written by C.~McMullen for the ICM Proceedings.

\section*{Brief biographical note}

{\small  Maryam Mirzakhani was born in 1977 in Teheran. After passing
an  entrance  test  she was accepted for the Farzanegan middle school
and  then  high  school  for  girls in Teheran administered by Iran's
National  Organization  for  Development  of  Exceptional Talents. As
Maryam   told  to  Erica~Klarreich,  journalist  of  Quanta  Magazine
(see~\cite{Quanta}  for  details)  she,  and her friend Roya Beheshti
convinced  the  principal of their high school to organize classes of
preparation  for  International  Mathematical  Olympiad (at this time
Iranian  team  never  contained  girls).  Maryam's determination gave
excellent  results: she won gold medals at International Mathematical
Olympiads in 1994 and in 1995.

Maryam   Mirzakhani   completed   undergraduate   studies  at  Sharif
University  in  Tehran  in  1999  and  moved  for graduate studies to
Harvard  University  choosing  Curtis McMullen as scientific advisor.
The   Ph.D.   Thesis  defended  in  2004  brought  her  international
recognition.  In her Thesis Maryam proved that the number $N(X,L)$ of
simple   closed   geodesics   (i.e.  of  those,  which  do  not  have
self-intersections)  on  a  Riemann surface $X$ of length at most $L$
grows asymptotically as
$$
N(X,L)\sim const(X)\cdot L^{6g-6} \text{ as } L\to\infty\,.
$$
(It  is  a  classical  result  that the number of \textit{all} closed
geodesics  grows  much  faster,  as $e^L/L$: most of closed geodesics
intersect  themselves.)  Another  extraordinary  result  of  the same
Ph.~D.  Thesis  was a new proof of Witten conjecture (the first proof
was obtained by M.~Kontsevich in 1992).

Having  defended the Ph.D.~Thesis Maryam Mirzakhani got a prestigious
Clay  Mathematics  Institute  Research  Fellowship  which  provides a
generous  salary  and  research expenses leaving to a Fellow complete
freedom  of choice where to perform research\footnote{Note that three
out  of  four  2014  Fields  Medalists  are also former Clay Research
Fellows.}.

Maryam~Mirzakhani  worked  for several years at Princeton University;
in  2008,  at  the age of 31, she became a Full Professor of Stanford
University, where she has been working ever since. Maryam is a mother
of a charming 3 years old daughter.

Personally  Maryam  is extremely nice and friendly, and not the least
bit  standoffish;  meeting  her at a conference you would take her at
the  first  glance  for  a young postdoc rather than for a celebrated
star.  All  her lectures which I have attended were full of sparkling
and  contagious  enthusiasm,  optimism, and appreciation of beauty of
mathematics;  they  inspire  you  to  attack  fearlessly  complicated
problems and, following Maryam, not to give up when they resist.
}

\section{From problems of solid state physics to surface foliations}
\label{s:phisics:2:foliations}

In  a sense, dynamical systems concern anything which moves; usually,
when  the  motion  has  already  achieved some kind of stable regime.
``The  thing  that  moves'' might be the Solar system, or a system of
particles  in  a  chamber,  or  a billiard ball on a table (where the
table  is  not  necessarily  a  rectangular  one), or currents in the
ocean,  or  electrons in a metal, etc. One can observe certain common
phenomena in large classes of dynamical systems; in particular, ideal
billiards  might  be interpreted as toy models of a gas in a chamber.
Such  toy models allow to elaborate tools to study original dynamical
systems of physical nature.

The   language  of  measured  foliations  on  surfaces  (generalizing
irrational  winding  lines  on  a  torus)  developed by Bill Thurston
proved  to  be  very  useful  in  working with the class of dynamical
systems   including   periodic  billiards  in  the  plane  (like  the
\textit{``windtree model''} introduced by Paule and Tatiana Ehrenfest
a  century  ago;  see the paper of P.~Hubert and R.~Krikorian in this
issue  for details) and dynamics of an electron on a Fermi-surface in
the   presence   of   a  homogeneous  magnetic  field  (\textit{``the
S.~P.~Novikov  Problem''};  see the survey~\cite{Maltsev:Novikov} for
details).

\begin{figure}[hbt]
\centering
\includegraphics{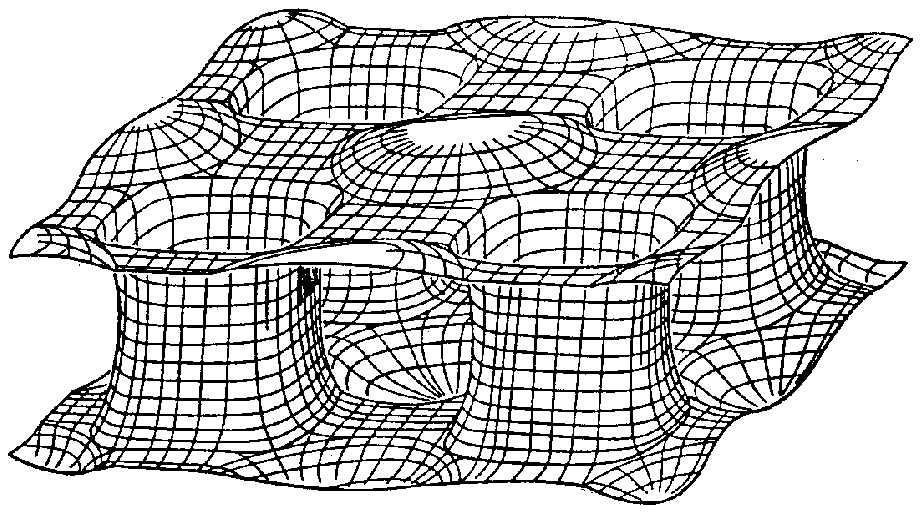}
\includegraphics{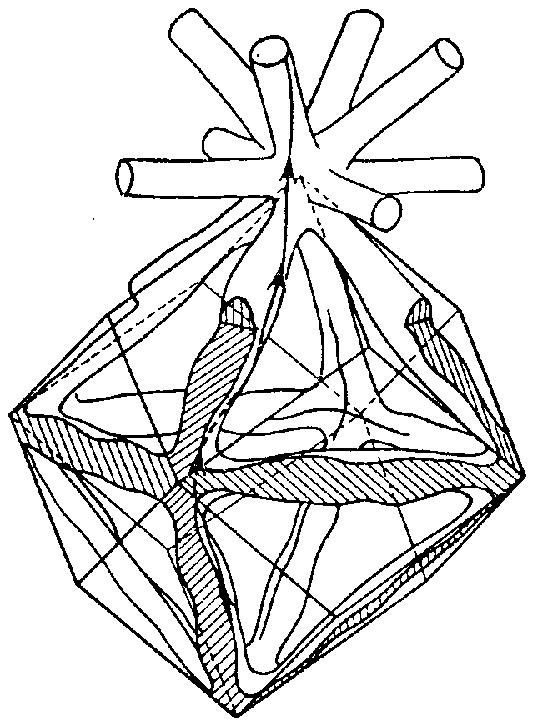}
\includegraphics{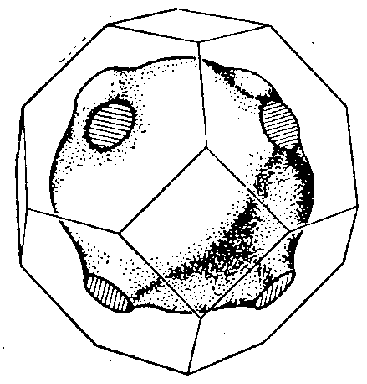}
\vspace{95bp} 
\caption{
\label{zorich:fig:fermi:surfaces}
Fundamental   domains  of  Fermi  surfaces  of  tin,  iron  and  gold
reproduced  from~\cite{zorich:LAK}.  Electron  trajectories  in these
metals  in the presence of a homogeneous magnetic field correspond to
plane sections of the corresponding periodic surfaces.
}
\end{figure}

A  flow following  the leaves of such measured foliation on a surface
decomposes the surface into two types of domains: periodic components
filled  with  periodic  trajectories  and minimal components in which
every  trajectory  is  dense.  After  applying a natural surgery, the
foliation on any minimal component can be in certain sense (described
in  the  next  section)  ``globally  straightened'':  one  can choose
appropriate flat coordinates in which the foliation is represented by
a   family   of   vertical   lines   $x=const$.   Versions   of  this
``straightening  theorem''  were  proved in different contexts and in
different        terms       by       E.~Calabi~\cite{zorich:Calabi},
A.~Katok~\cite{zorich:Katok:1973},                  \mbox{J.~Hubbard}
and~H.~Masur~\cite{zorich:Hubbard:Masur}, and by other authors. These
theorems  imply that to study an important class of dynamical systems
it  is  sufficient  to confine the study to a simplified model and to
study straight line foliations on ``translation surfaces''.

\section{Translation surfaces}
\label{s:Translation:surfaces}

In  this  section  we  describe  a  very  concrete  construction of a
vertical  foliation  on  a surface endowed with a rather special flat
metric  with  isolated  singularities.  The ``straightening theorem''
mentioned  above  asserts  that all orientable measured foliations on
surfaces  can  be reduced (in the way described above) to the ones as
in the current section.

\begin{figure}[htb]
%
%
\includegraphics{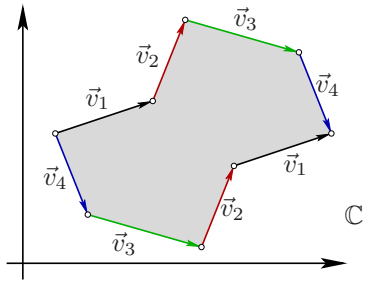}
\begin{picture}(35,2)(142,2) 
\put(122,-35){$\vec v_1$}
\put(141,-19){$\vec v_2$}
\put(178,-5){$\vec v_3$}
\put(209,-30){$\vec v_4$}
\put(106,-65){$\vec v_4$}
\put(132,-90){$\vec v_3$}
\put(173,-77){$\vec v_2$}
\put(197,-60){$\vec v_1$}
\put(220,-80){$\C{}$}
\end{picture}
\vspace{100bp}
\caption{
\label{fig:polygon}
Identifying  corresponding pairs of sides of this polygon by parallel
translations  we  obtain  a  flat surface of genus two endowed with a
flat metric having a single conical singularity.
}
\end{figure}

Consider  a  collection of vectors $ \vec v_1, \dots, \vec v_\noi$ in
$\R{2}$  and  arrange  these  vectors  into  a broken line. Construct
another  broken  line  starting  at  the  same point as the first one
arranging the same vectors in the order $\vec v_{\pi(1)}, \dots, \vec
v_{\pi(\noi)}$,  where  $\pi$ is some permutation of $\noi$ elements.
By  construction  the  two  broken  lines  share  the same endpoints;
suppose  that  they  bound  a polygon as in Figure~\ref{fig:polygon}.
Identifying  the  pairs  of  sides  corresponding to the same vectors
$\vec  v_j$, $j=1, \dots, \noi$, by parallel translations we obtain a
closed topological surface.

By  construction,  the  surface  is  endowed with a flat metric. When
$\noi  =  2$  and  $\pi=(2,1)$ we get a usual flat torus glued from a
parallelogram.  For  larger number of elements we might get a surface
of  higher  genus,  where  the genus is determined by the permutation
$\pi$.   It   is  convenient  to  impose  from  now  on  some  simple
restrictions on the permutation $\pi$ which guarantee, in particular,
non  degeneracy of the surface; see~\cite{Masur} or~\cite{Veech}.

\begin{figure}[hbt]
\centering
\includegraphics{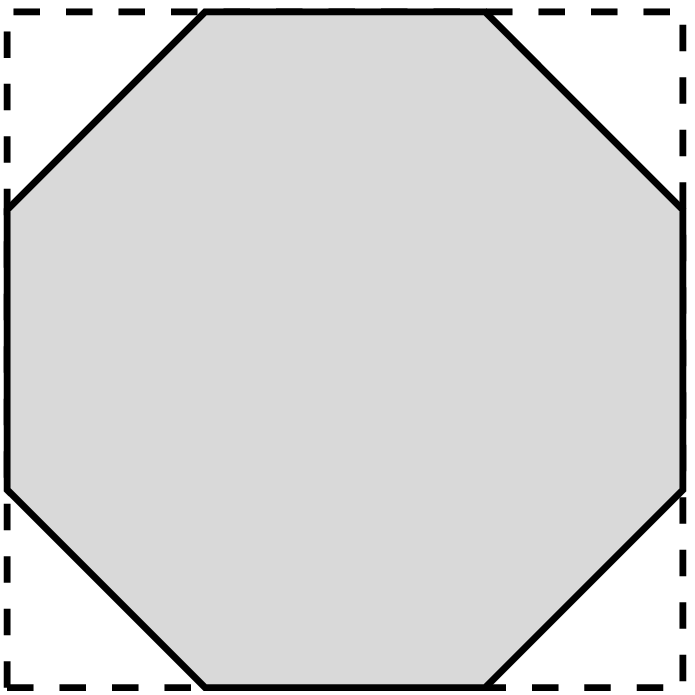}
\includegraphics{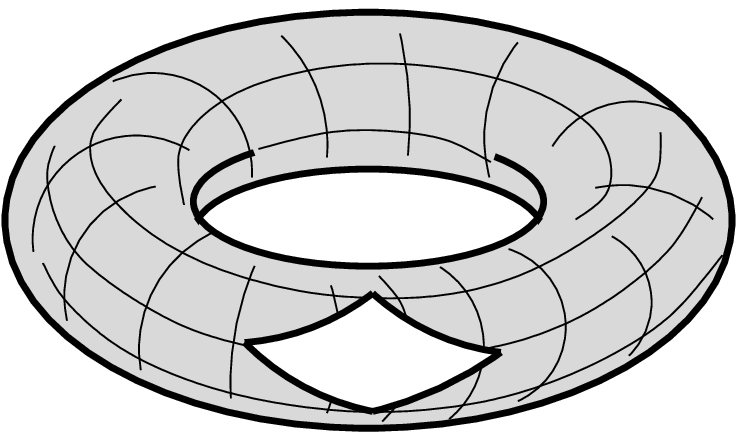}
\includegraphics{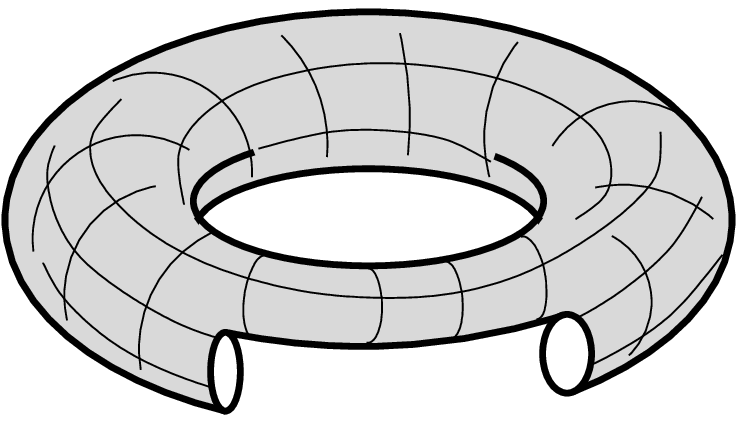}
\includegraphics{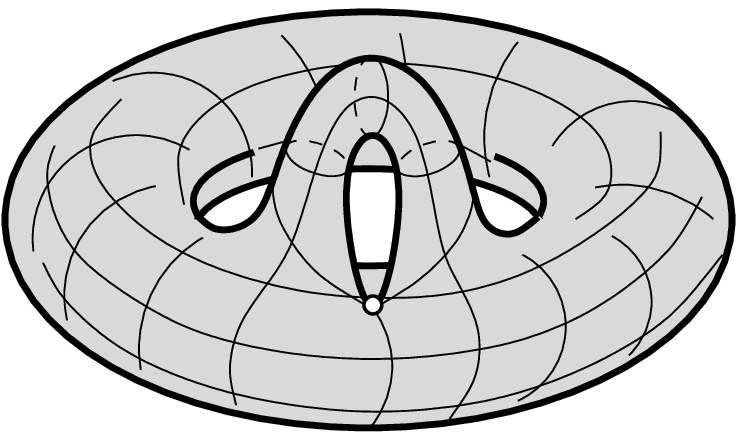}
\vspace{50bp}
\caption{
\label{fig:cartoon}
Cartoon movie of gluing a translation surface of genus
two  from a regular octagon.
}
\end{figure}

For  example,  a regular octagon gives rise to a surface of genus two
as   in   Figure~\ref{fig:cartoon}.   Indeed,  identifying  pairs  of
horizontal  and  vertical  sides  of a regular octagon we get a usual
torus  with  a hole in the form of a square. We slightly cheat in the
next  frame, where we turn this hole by $45^\circ$ and only then glue
the  next pair of sides. As a result we get a torus with two isolated
holes  as on the third frame. Identifying the remaining pair of sides
(which  represent  the  holes) we get a torus \textit{with a handle},
or, in other words, a surface of genus two.

Similar to the torus case, the surface glued from the regular octagon
or  from an octagon as in Figure~\ref{fig:polygon} also inherits from
the  polygon  a  flat metric, but now the resulting flat metric has a
singularity  at  the  point  obtained from identified vertices of the
octagon.

Note  that the flat metric thus constructed is very special: since we
identify  the sides of the polygon only by translations, the parallel
transport  of  any  tangent  vector  along  a  closed cycle (avoiding
conical  singularities)  on  the  resulting surface brings the vector
back  to  itself. In other words, our flat metric has \textit{trivial
holonomy}.  In  particular,  since a parallel transport along a small
loop  around any conical singularity brings the vector to itself, the
cone  angle  at  any singularity is an integer multiple of $2\pi$. In
the  most  general situation the flat surface of genus $g$ would have
several  conical  singularities  with  cone angles $2\pi(d_1+1),\dots
2\pi(d_\noz+1)$, where $d_1+\dots +d_\noz=2g-2$.

It  is  convenient  to consider the vertical direction as part of the
structure. A surface endowed with a flat metric with trivial holonomy
and   with   a   choice   of   a   vertical  direction  is  called  a
\textit{translation  surface}. Two polygons in the plane obtained one
from  another  by  a  parallel  translation  give  rise  to  the same
translation  surface,  while  polygons obtained one from another by a
nontrivial  rotation  (usually)  give  rise  to  distinct translation
surfaces.

We  can  assume  that the polygon defining our translation surface is
embedded  into  the  complex  plane $\C{}\simeq\R{2}$ with coordinate
$z$.   The   translation   surface   obtained   by   identifying  the
corresponding  sides  of  the polygon inherits the complex structure.
Moreover,  since  the  gluing  rule for the sides can be expressed in
local  coordinates  as  $z=\tilde z+const$, the closed 1-form $dz$ is
well-defined not only in the polygon, but on the surface. An exercise
in  complex  analysis shows that the complex structure extends to the
points  coming  from the vertices of the polygon, and that the 1-form
$\omega=dz$  extends  to  the  holomorphic  1-form  on  the resulting
Riemann  surface.  This  1-form  $\omega$ has zeroes of degrees $d_1,
\dots,  d_\noz$  exactly  at  the  points  where  the flat metric has
conical singularities of angles $2\pi(d_1+1), \dots, 2\pi(d_\noz+1)$.

Reciprocally,  given  a  holomorphic  1-form  $\omega$  on  a Riemann
surface   one   can   always  find  a  local  coordinate  $z$  (in  a
simply-connected  domain not containing zeroes of $\omega$) such that
$\omega=dz$.  This  coordinate is defined up to an additive constant.
It  defines  the translation structure on the surface. Cutting up the
surface  along an appropriate collection of straight segments joining
conical  singularities  we  can  unwrap  the  Riemann  surface into a
polygon as above.

This  construction  shows  that  the  two  structures  are completely
equivalent;  the  flat  metric with trivial holonomy plus a choice of
distinguished  direction or a pair: Riemann surface and a holomorphic
1-form on it.

\section{Families of translation surfaces and dynamics
in the moduli space}
\label{s:modulis:space}

The  polygon  in our construction depends continuously on the vectors
$\vec{v}_i$.   This   means   that  the  topology  of  the  resulting
translation  surface  (its genus $g$, the number and the types of the
resulting   conical   singularities)  does  not  change  under  small
deformations of the vectors $\vec{v}_i$. For every collection of cone
angles     $2\pi(d_1+1),     \dots,     2\pi(d_\noz+1)$    satisfying
$d_1+\dots+d_\noz=2g-2$  with integer $d_i$ for $i=1,\dots, \noi$, we
get   a   \textit{family   $\cH(d_1,\dots,d_\noz)$   of   translation
surfaces}.  Vectors  $\vec  v_1, \dots, \vec v_\noi$ can be viewed as
complex  coordinates  in  this  space,  called  \textit{period
coordinates}.   These   coordinates   define   a   structure   of   a
\textit{complex  orbifold}  (manifold with moderate singularities) on
each  space  $\cH(d_1,\dots,d_\noz)$.  The  geometry  and topology of
spaces of translation surfaces is not yet sufficiently explored.

Readers  preferring  algebro-geometric  language may view a family of
translation  surfaces  with fixed conical singularities $2\pi(d_1+1),
\dots,  2\pi(d_\noz+1)$ as the stratum $\cH(d_1,\dots,d_\noz)$ in the
moduli  space  $\cH_g$  of  pairs  (Riemann  surface $C$; holomorphic
1-form  $\omega$  on  $C$),  where  the  stratum  is specified by the
degrees   $d_1,\dots,   d_\noz$   of   zeroes   of   $\omega$,  where
$d_1+\dots+d_\noz=2g-2$.  Note that while the moduli space $\cH_g$ is
a holomorphic $\C{g}$-bundle over the moduli space $\cM_g$ of Riemann
surfaces, individual strata are not. For example, the minimal stratum
$\cH(2g-2)$  has  complex  dimension  $2g$,  while  the  moduli space
$\cM_g$  has  complex  dimension  $3g-3$.  The  very  existence  of a
holomorphic  form  with  a  single zero of degree $2g-2$ on a Riemann
surface $C$ is a strong condition on $C$.

To  complete  the description of the space of translation surfaces we
need  to present one more very important structure: the action of the
group  $\GL$  on  $\cH_g$  preserving strata. The description of this
action  is  particularly simple in terms of our polygonal model $\Pi$
of  a  translation  surface $S$. A linear transformation $g\in\GL$ of
the plane maps the polygon $\Pi$ to a polygon $g\Pi$. The new polygon
again  has all sides arranged into pairs, where the two sides in each
pair  are  parallel  and  have  equal  length.  We  can  glue  a  new
translation  surface  and  call it $g\cdot S$. It is easy to see that
unwrapping  the  initial  surface  into  different polygons would not
affect the construction. Note also, that we explicitly use the choice
of  the  vertical direction: any polygon is endowed with an embedding
into $\R{2}$ defined up to a parallel translation.

The  subgroup  $\SL\subset\GL$ preserves the flat area. This implies,
that   the   action   of   $\SL$   preserves  the  real  hypersurface
$\cH_1(d_1,\dots,d_\noz)$  of translation surfaces of area one in any
stratum   $\cH(d_1,\dots,d_\noz)$.   The   codimension-one   subspace
$\cH_1(d_1,\dots,d_\noz)$  can  be  compared  to  the unit sphere (or
rather   to   the   unit   hyperboloid)   in   the   ambient  stratum
$\cH(d_1,\dots,d_\noz)$.

Recall  that  under  appropriate assumption on the permutation $\pi$,
the $\noi$ vectors
$$
\vec v_1=\begin{pmatrix}v_{1,x}\\v_{1,y}\end{pmatrix}       \dots,
\vec v_\noi=\begin{pmatrix}v_{\noi,x}\\v_{\noi,y}\end{pmatrix}
$$
as  in  Figure~\ref{fig:polygon}  define  local  coordinates  in  the
embodying family $\cH(d_1,\dots,d_\noz)$ of translation surfaces. Let
$d\nu:=dv_{1   x}dv_{1  y}\dots  dv_{\noi  x}d  v_{\noi  y}$  be  the
associated  volume  element  in  the  corresponding  coordinate chart
$U\subset\R{2\noi}$.  It  is  easy  to  verify,  that $d\nu$ does not
depend on the choice of ``coordinates'' $\vec v_1,\dots,\vec v_\noi$,
so  it  is  well-defined on $\cH(d_1,\dots,d_\noz)$. Similarly to the
case  of  Euclidian  volume  element, we get a natural induced volume
element  $d\nu_1$  on the unit hyperboloid $\cH_1(d_1,\dots,d_\noz)$.
It  is easy to check that the action of the group $\SL$ preserves the
volume element $d\nu_1$.

The  following  Theorem  proved  independently  and simultaneously by
\mbox{H.~Masur}~\cite{Masur}   and   W.~Veech~\cite{Veech}   is   the
keystone of this area.
\begin{Theorem}[H.~Masur; W.~A.~Veech]
The  total  volume  of  every  stratum  $\cH_1(d_1,\dots,d_\noz)$  is
finite.

The  group  $\SL$  and its diagonal subgroup act ergodically on every
connected component of every stratum $\cH_1(d_1,\dots,d_\noz)$.
\end{Theorem}

\begin{figure}[htb]
%
\includegraphics{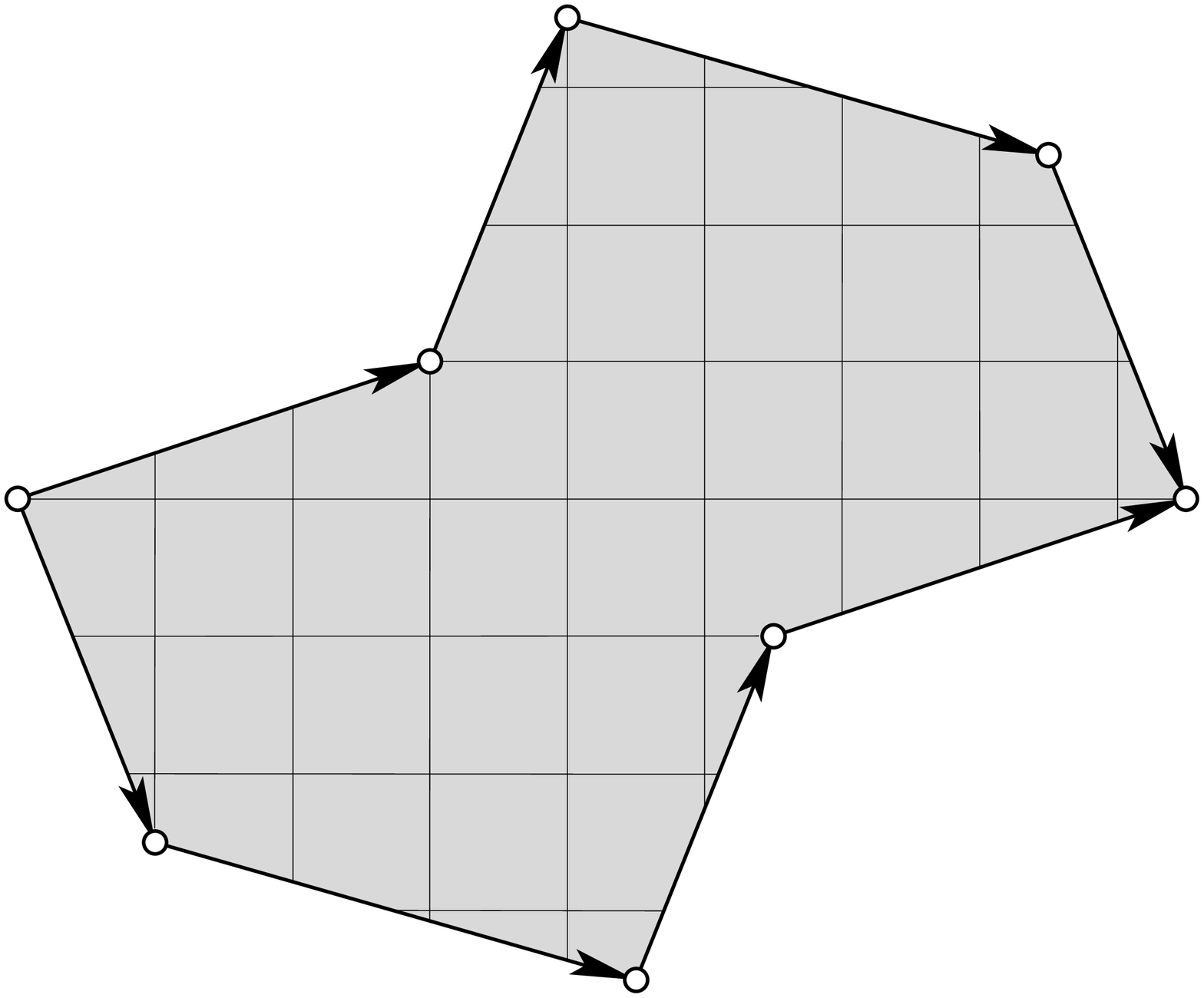}
\includegraphics{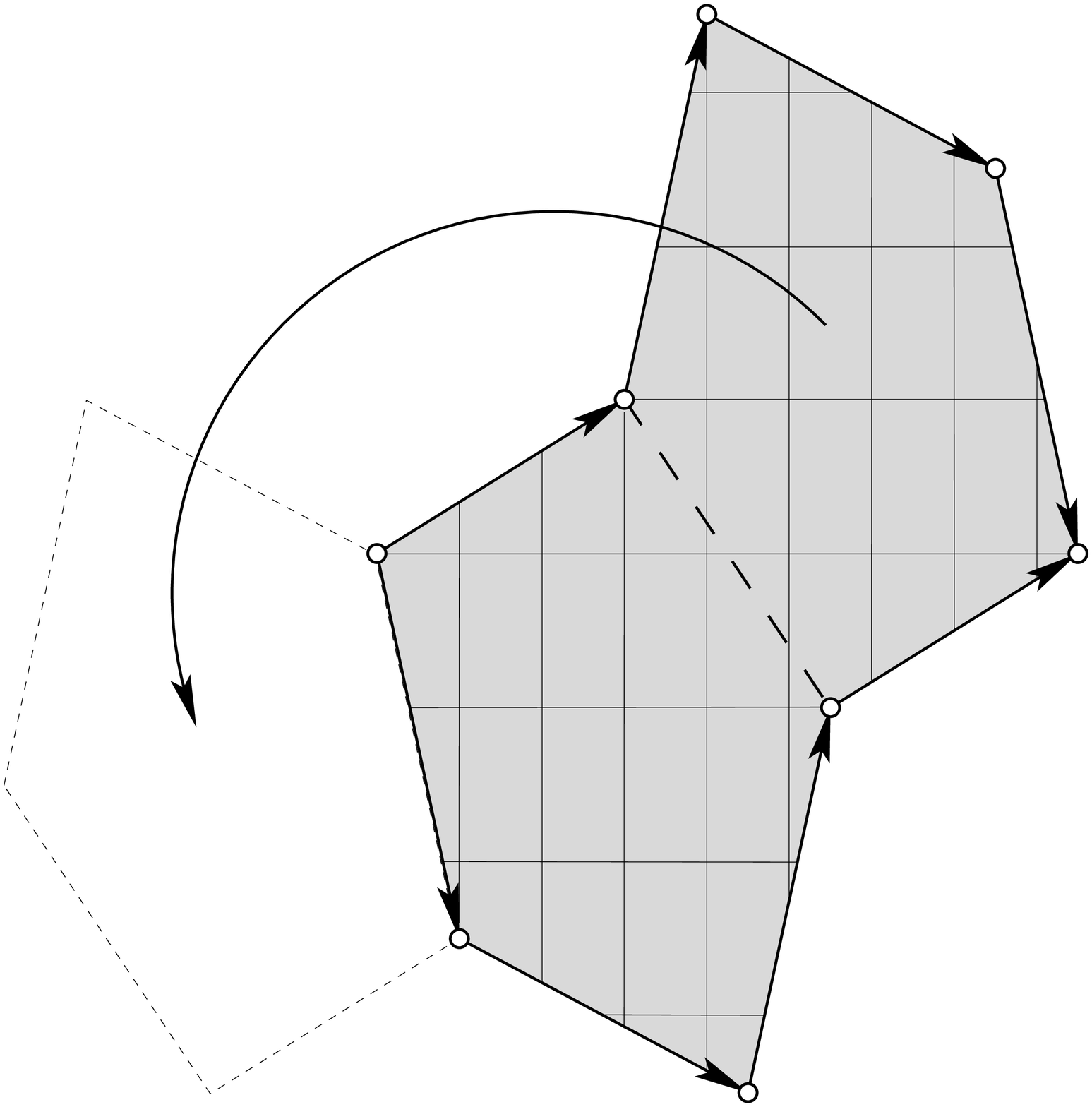}
\includegraphics{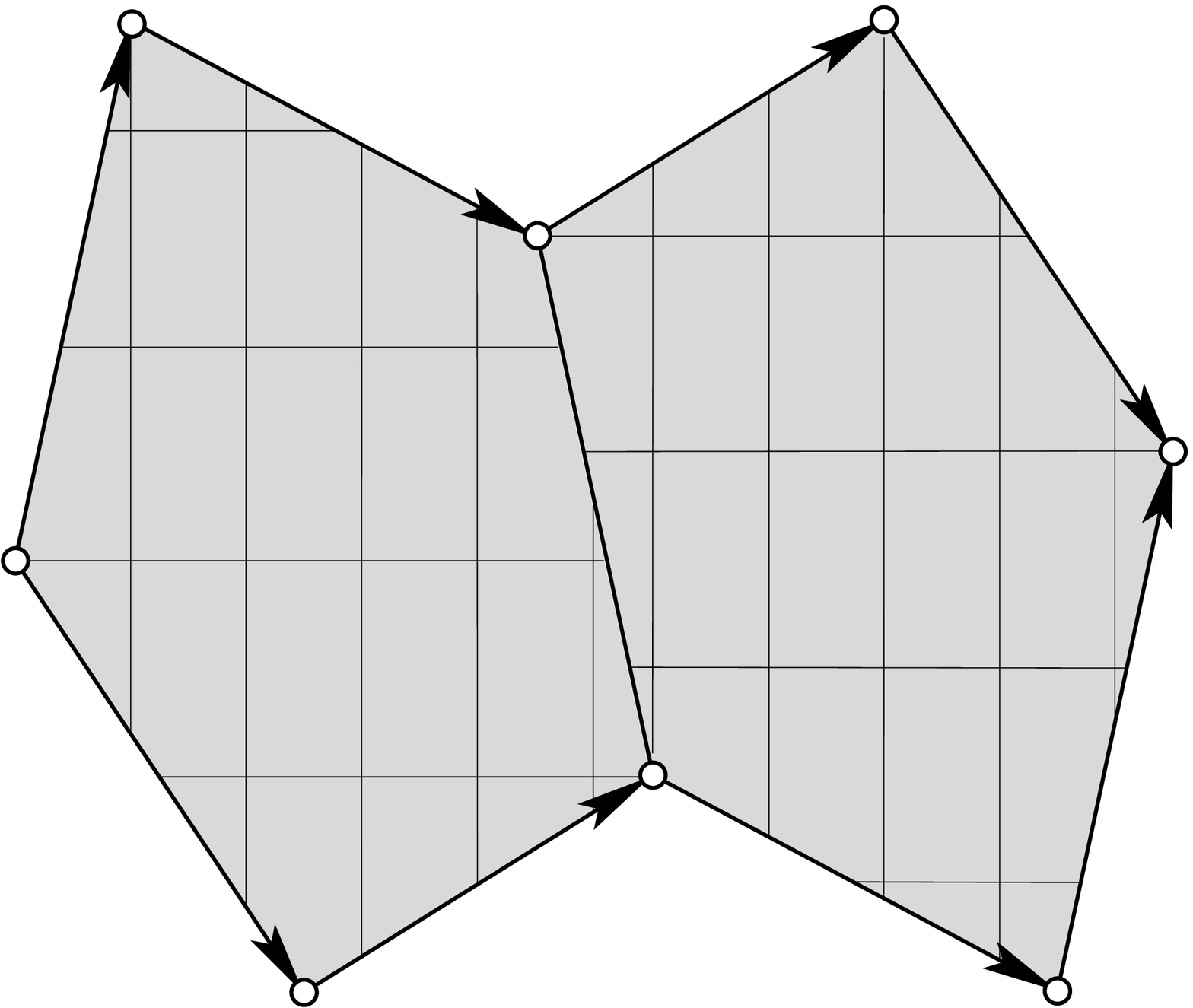}
\includegraphics{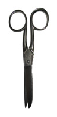}
\begin{picture}(0,0)(-5,0)
\put(-95,-70){\LARGE$\longrightarrow$}
\put(40,-70){\LARGE$\ =$}
\end{picture}
\vspace{100bp}
\caption{
\label{fig:teich:flow}
Note  that  expansion-contraction  (action  of  the  diagonal  group)
changes  the translation surface, white cut-and-paste transformations
change  only  the  polygonal  pattern,  and  do  not  change the flat
surface.
}
\end{figure}

Here ``ergodically'' means that any measurable subset invariant under
the action of the group has necessarily measure zero or full measure.
Ergodic  theorem  claims  that  in such situation the orbit of almost
every  point  homogeneously fills the ambient connected component. In
plain  terms,  the  ergodicity of the action of the diagonal subgroup
can be interpreted as follows. Having almost any polygon as above, we
can  choose appropriate sequence of times $t_i$ such that contracting
the  polygon  horizontally  with  a factor $e^{t_i}$ and expanding it
vertically with the same factor $e^{t_i}$ and modifying the resulting
polygonal   pattern  of  the  resulting  translation  surface  by  an
appropriate    sequence   of   cut-and-paste   transformations   (see
Figure~\ref{fig:teich:flow})  we  can  get  arbitrary  close to, say,
regular octagon rotated by any angle chosen in advance.

Now everything is prepared to present the first marvel of this story.
Suppose that we want to find out some fine properties of the vertical
flow  on  an  individual translation surface. Applying a homothety we
can  assume  that  the translation surface has area one. Howard~Masur
and  William~Veech  suggest  the  following  approach.  Consider  our
translation  surface  (endowed  with a vertical direction) as a point
$S$  in  the  ambient stratum $\cH_1(d_1,\dots,d_\noz)$. Consider the
orbit  $\SL\cdot  S$  (or  the  orbit  of $S$ under the action of the
diagonal    subgroup    $\begin{pmatrix}e^t&0\\0&e^{-t}\end{pmatrix}$
depending  on  the initial problem). Numerous important properties of
the  initial  vertical flow are encoded in the geometrical properties
of  the  closure of the corresponding orbit. This approach placed the
problem  of finding the orbit closures under the action of $\SL$, and
studying their geometry at the center of the studies in this area for
the last three decades. To give at least one example of this approach
we  present  Masur's  criterion  of \textit{unique ergodicity} of the
vertical flow.

\begin{figure}[hbt]
   %
\includegraphics{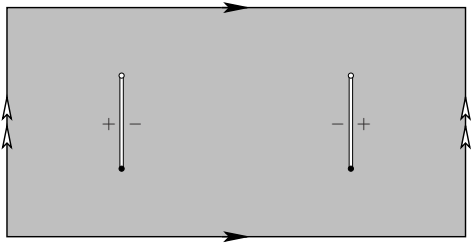}
\begin{picture}(0,0)(0,0)
\put(-50,-49){\tiny $+$}
\put(-40,-49){\tiny $-$}
\put(36.5,-49){\tiny $-$}
\put(46.5,-49){\tiny $+$}
\end{picture}
\vspace{90bp}
\caption{
\label{fig:nonergodic}
For  appropriate  directions  of  the  straight  line  flow  on  this
translation  surface  all trajectories are dense, but the flow is not
ergodic:  trajectories  are \textit{not} uniformly distributed on the
translation surface.
}
\end{figure}

In  the  late  60's  and  70's  several  authors (including W.~Veech,
M.~Keane,  \mbox{A.~Katok}, and others) discovered in different terms
and  in  various  contexts the phenomenon which can be illustrated by
the   following   example.   Take   a   translation   surface  as  in
Figure~\ref{fig:nonergodic},  where  the  rectangle has size $1\times
2$,  and the slits have any irrational length. (We are not restricted
to  polygons  as in Figure~\ref{fig:polygon} to construct translation
surfaces;  we  could  also  glue  translation surfaces from triangles
imbedded  into the Eucledian plane, provided that all identifications
of  sides  are  parallel  translations.)  For  uncountable  number of
directions all trajectories of the straight line flow would be dense,
but  it  would  not be ergodic: some trajectories would spend most of
the  time  in  the  middle part of the translation surface, while the
other     ones    ---    mostly    in    the    complementary    part
(see~\cite{Masur:Handbook:1B}  for  details).  Clearly,  applying  an
appropriate  rotation  to  the  translation  surface,  one can make a
nonergodic direction vertical.

Having  presented  the  phenomenon, we can state Masur's criterion of
unique ergodicity of the vertical flow.

\begin{Theorem}[H.~Masur]
Suppose  $S$  is  a  translation  surface.  Suppose  the  flow in the
vertical  direction  on  $S$  is  not  uniquely ergodic. Then for any
compact  set  $K\subset  \cH(d_1,...,d_\noz)$ in the ambient space of
translation  surfaces  the trajectory $g_t\cdot S$ of $S$ would never
visit $K$ for $t$ sufficiently large, $t>t_0(K)$. In other words, the
trajectory   of  $S$  under  the  action  of  the  diagonal  subgroup
$\begin{pmatrix}e^t&0\\0&e^{-t}\end{pmatrix}$  eventually  leaves any
compact   set;   it   ``escapes   to  a  multidimensional  cusp''  of
$\cH(d_1,...,d_\noz)$.

\noindent
(Actually,  the  statement  is  even  stronger: the projection of the
trajectory  $g_t\cdot  S$  to  the  moduli  space  $\cM_g$ of Riemann
surfaces eventually leaves any compact set in $\cM_g$.)
\end{Theorem}

The  vertical  flow  on the surface in Figure~\ref{fig:nonergodic} is
periodic,  so  extremely  nonergodic.  And  indeed, the action of the
diagonal  subgroup  makes  our  flat torus with slits become more and
more narrow and long. In this case the family of deformed translation
surfaces  leaves  any  compact  set  in a very simple way: there is a
closed  flat  geodesics  which  gets  pinched.  To  construct minimal
nonergodic  examples mentioned above one should show that for certain
angles the vertical flow applied to rotated translation surface makes
the  deformed translation surface leave eventually any compact set in
the family of flat surfaces and never return back.

At  the  first  glance,  we  have  just reduced the study of a rather
simple   dynamical   system,  namely,  of  the  vertical  flow  on  a
translation  surface, to a really complicated one --- to the study of
the  action  of the group $\GL$ on the space $\cH(d_1,\dots,d_\noz)$.
Nevertheless,  this  approach  proved  to  be  extremely fruitful and
powerful,  despite  the  fact that geometry and topology of spaces of
translation surfaces is not sufficiently explored yet.

However,  as  it  would  be explained in the next section, successful
implementations of this approach were often based on a chain of happy
coincidences (like reduction to translation surfaces of genera one or
two  which  are very special). The Magic Wand Theorem of A.~Eskin and
M.~Mirzakhani and of A.~Eskin, M.~Mirzakhani, A.~Mohammadi transforms
this approach from art to science.

\section{``Almost all'' versus ``all''}
\label{s:almost:all}

Even  for  those  classes of dynamical systems which are sufficiently
well  understood, the only kind of predictions of ``what would happen
to  a  particle  after  sufficiently  long time'' always contain some
version  of  the  word  ``typically''  usually  meaning  ``for a full
measure  set  of initial data''. The trouble (which, depending on the
taste,  might be considered as an advantage: do not get distracted by
details) is that even for those dynamical systems which are very well
studied  and  understood,  and where one knows, basically, everything
about  ``typical  behavior''  of  trajectories,  one  can  say almost
nothing  about  behavior of any concrete particular trajectory: there
is  no  way  to  tell,  whether  your  particular  starting  data are
``typical''  or  not.  If  you  repeat  thousands of experiments with
random  starting  data and you want to establish some statistics, you
do  not  care  about  rare  nontypical  fluctuations.  But if you are
interested  in  the  future  of some very special asteroid B 612, and
only  by  this,  most  of  the  methods  of  dynamical systems become
completely useless for you.

The  difficulty  is  conceptual;  it  is  neither  related to lack of
knowledge  at the current state of development of mathematics, nor to
the  presence  of  noise,  or  friction,  etc  in realistic dynamical
systems. Even for absolutely deterministic systems, and even assuming
all  necessary mathematical abstractions like absence of any noise or
friction,  the  trouble  persists.  The  reason  is that for the vast
majority  of  dynamical  systems  (in particular, for very smooth and
nice   ones)  certain  individual  trajectories  might  be  extremely
sophisticated.  For example, they can cover extremely fractal sets on
a large scale of time.

For  example,  the map $f:x\mapsto \{2x\}$ homogeneously twisting the
unit  circle  $S^1=\R{}/\Z{}$  twice around itself has orbits filling
Cantor sets of, basically, arbitrary Hausdorff dimension between zero
and  one;  has  nonclosed orbits avoiding certain arcs of the circle,
etc.   In   other  words,  this  extremely  nice  map,  clearly,  has
trajectories with very peculiar properties.

All  these  properties  become  much  more  visible  using the binary
representation  of  a  real  number  $x\in[0,1[$ instead of the usual
decimal one. If
$$
x=\frac{n_1}{2}+\dots+\frac{n_k}{2^k}+\cdots\,,
$$
where  all  binary  digits $n_k$ are zeroes or ones, then the map $f$
acts on the sequence $(n_1,n_2,\dots,n_k,\dots)$ by erasing the first
digit.  (This  operation  on  the space of semi-infinite sequences of
zeroes and ones is called the \textit{Bernoulli shift}).

The geodesic flow on any compact Riemann surface of constant negative
curvature   has  similar  behavior.  It  was  observed  long  ago  by
H.~Furstenberg   and   B.~Weiss   that  the  closures  of  individual
trajectories  might  have  arbitrary  (or almost arbitrary) Hausdorff
dimension  in  the  range  from 1 (closed trajectories) to 3 (typical
trajectories).

A  straight-line flow on a torus $\T{n}=\R{n}/\Z{n}$ is an example of
a  (very  rare)  dynamical  system, where the closure of \textit{any}
orbit is a nice submanifold. Let $\vec{V}=(V_1,\dots,V_n)$ denote the
direction  of  the  flow.  The  closure of \textit{any} trajectory in
direction  $\vec{V}$  is  a sub-torus $\T{k}$, where $1\le k\le n$ is
the  degree  of irrationality $k=\dim_{\Q{}}\{V_1,\dots,V_n\}$ of the
direction.  Say,  in  the particular case of a two-dimensional torus,
when  $n=2$, all trajectories of the flow in a rational direction are
already  closed  --- they are circles $S^1=\T{1}$, and the closure of
any  trajectory  of the flow in an irrational direction is the entire
torus  $\T{2}$.  Actually,  this is not surprising at all: torus is a
homogeneous  space,  and  the  group  of  automorphisms  of the torus
preserving the flow acts transitively on the torus.

In  a  sense, up to know, there was only one known class of dynamical
systems,  for  which  one  could  find  the  closure  of  any  single
trajectory,  and  for which all possible closures were described by a
short  list  of  possible  simple cases like in the example above. It
happens for very special dynamical systems in homogeneous spaces. One
of  the  key  statements  in this theory was proved by Marina~Ratner;
extremely important contributions to this theory as well as fantastic
applications  to  the  number  theory,  were developed by S.~G.~Dani,
G.~Margulis,  and  by other great mathematicians, including A.~Eskin,
S.~Mozes,  and  N.~Shah.  The scale of applications of this theory to
different   areas   of   mathematics  continues  to  extend.  Indeed,
homogeneous   spaces   naturally   appear   in   various  domains  of
mathematics.  (Both  the  theory  and  the list of major contributors
merit a separate paper rather than a short paragraph.)

\section{Magic Wand Theorem}
\label{s:Theorem}

And    now    comes   the   result   of   Alex   Eskin   and   Maryam
Mirzakhani~\cite{Eskin:Mirzakhani}  (incorporating  the joint results
of             these             authors            and            of
Amir~Mohammadi~\cite{Eskin:Mirzakhani:Mohammadi}).

It  is  known  that  the  moduli  space \textit{is not} a homogeneous
space.  Nevertheless,  the  orbit  closures  of $\GL$ in the space of
translation  surfaces  are  as  nice  as  one can only hope: they are
complex   manifolds   possibly   with   very  moderate  singularities
(so-called  ``orbifolds''). In this sense the action of the $\GL$ and
of  $\SL$ on the space of translation surfaces described above mimics
certain  properties  of  the  dynamical systems in homogeneous spaces
mentioned at the end of the previous section.

\begin{MagicWand}
The  closure  of  any  $\GL$-orbit is a complex suborbifold (possibly
with   self-intersections);   in   period   coordinates  $\vec
v_1,\dots,     \vec    v_\noi$    in    the    corresponding    space
$\cH(d_1,\dots,d_\noz)$   of   translation   surfaces  it  is  locally
represented by an affine subspace.

Any ergodic $\SL$-invariant measure is supported on a suborbifold. In
coordinates  $\vec  v_1,  \dots,  \vec  v_\noi$  this  suborbifold is
represented  by an affine subspace, and the invariant measure is just
a usual affine measure on this affine subspace.
\end{MagicWand}

As  a  vague  conjecture (or, better say, as a very optimistic dream)
this property was discussed since long ago, and since long ago, there
was  not  a slightest hint for a general proof. The only exception is
the  case  of  surfaces  of  genus  two,  for  which  ten  years  ago
C.~McMullen   proved   a   very   precise   statement~\cite{McMullen}
classifying  all possible orbit closures. He used, in particular, the
hard  artillery  of  Ratner's  results  which  are  applicable  here.
However,  the  theorem  of  McMullen  is  based  on  the very special
properties  of  surfaces  of  genus  two,  which do not generalize to
higher genera.

The proof of Alex Eskin and Maryam Mirzakhani is a titanic work which
took  many  years.  It  absorbed numerous fundamental developments in
dynamical  systems  which  do  not have any direct relation to moduli
spaces.  To  mention only a few, it incorporates certain ideas of low
entropy  method  of M.~Einsiedler, A.Katok, E.~Lindenstrauss; results
of  G.~Forni  and  of  M.~Kontsevich  on  Lyapunov  exponents  of the
Teichm\"uller  geodesic  flow;  the ideas from the works of Y.~Benoit
and  \mbox{J.-F.~Quint} on stationary measures; iterative improvement
of  the  properties of the invariant measure inspired by the approach
of  G.~Margulis  and  G.~Tomanov to the actions of unipotent flows on
homogeneous  spaces; some fine ergodic results due to Y.~Guivarch and
A.~Raugi.

What  can  we  do  now,  when this theorem is proved? For example, in
certain situations this theorem works like a Magic Wand, which allows
to  touch  \textit{any}  given billiard in certain class of billiards
and  in  theory  (and  more  and  more  often  in  practice) find the
corresponding  orbit  closure  in  the  moduli  space  of translation
surfaces.  The  geometry  of this orbit closure tells you, basically,
everything you want to know about the initial billiard.

\begin{figure}[hbt]
   %
\includegraphics{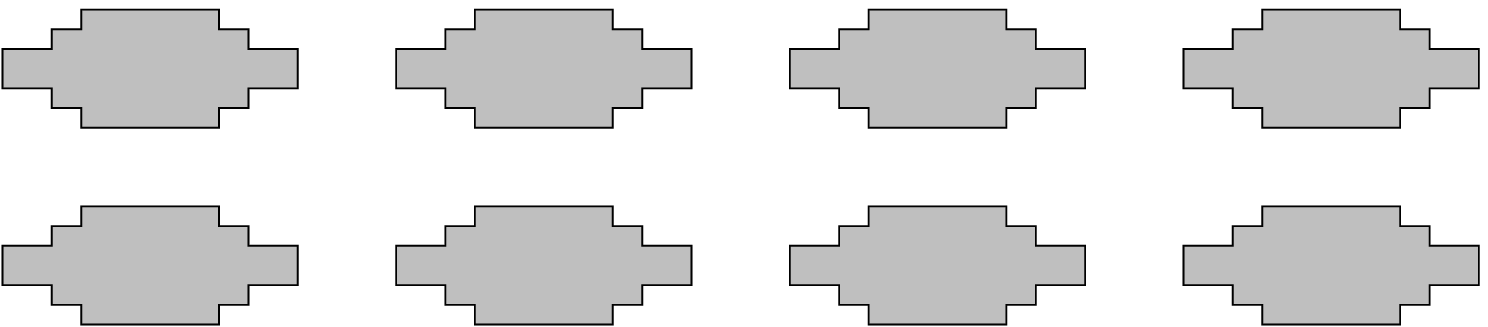}
\includegraphics{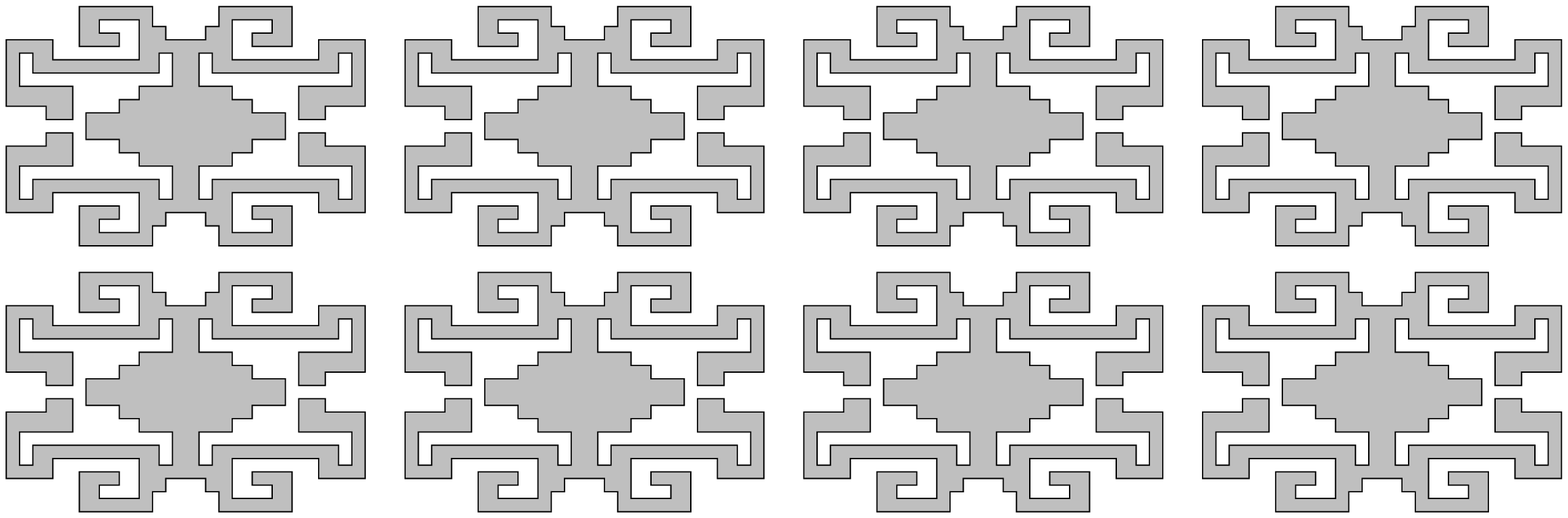}
\vspace{50bp}
\caption{
\label{fig:generalized:windtree}
The diffusion rate depends only on the number of the angles of almost
any symmetric obstacle as on the picture .
}
\end{figure}

Suppose   you  want  to  study  the  billiard  in  the  plane  filled
periodically        with        the       obstacles       as       in
Figure~\ref{fig:generalized:windtree} (see the paper of P.~Hubert and
R.~Krikorian~\cite{Hubert:Krikorian}  in this issue discussing such a
\textit{windtree model}). A trajectory might go far away, then return
relatively  close  back  to  the starting point, then make other long
trips.  The  \textit{diffusion rate} $\nu$ describes the average rate
$T^\nu$  with  which  the  trajectory  expands in the plane on a long
range of time $T\gg 1$. More formally,
$$
\nu:=\limsup_{T\to\infty}\frac{\log(\text{diameter  of  part of trajectory
for interval of time }[0,T])}{\log T}\,.
$$
For  the  usual  random  walk  in  the  plane, or for a billiard with
periodic  circular obstacles the diffusion rate is known to be $1/2$:
the  most  distant  point  of  a piece of trajectory corresponding to
segment  of  time  $[0,T]$  would  be  located  roughly at a distance
$\sqrt{T}$.  (We  do  not  discuss  at  what  time  $t_0\in[0,T]$ the
trajectory would be that far.)

The  Magic  Wand Theorem and companion results allow to show that for
\textit{any} obstacle as in Figure~\ref{fig:generalized:windtree} the
diffusion  rate $\nu$ is one and the same for all starting points and
for  almost  all  starting directions. Moreover, one can even compute
the  diffusion  rate!  To  perform  this  task  one has to proceed as
follows.  Find  the  associated translation surface $S$; it is really
easy.   Using   the  Magic  Wand  theorem  (and  its  development  by
A.~Wright~\cite{Wright})  find the corresponding orbit closure $\cL$.
By      a     very     recent     theorem     of     A.~Eskin     and
J.~Chaika~\cite{Eskin:Chaika}  almost all directions for \textit{any}
translation surface are Lyapunov-generic, so it is sufficient to find
the  appropriate Lyapunov exponent of the Teichm\"uller geodesic flow
on  $\cL$ and you are done! (See details on the windtree billiard and
on  Lyapunov  exponents  in the paper~\cite{Hubert:Krikorian} in this
issue.)

Say,  for  almost  any symmetric obstacle with $4m-4$ angles $3\pi/2$
and  $4m$  angles $\pi/2$ V.~Delecroix and the author showed that the
diffusion rate is
$$
\frac{(2m)!!}{(2m+1)!!}
\quad\sim\frac{\sqrt{\pi}}{2\sqrt{m}}\quad
\text{as }\ m\to\infty\,.
$$
(This  answers  a question by J.-C.~Yoccoz whether for certain shapes
of  the  obstacles  the  diffusion  rate  can  be arbitrary small and
develops     the     original     answer     for    $m=1$    obtained
in~\cite{Delecroix:Hubert:Lelievre}).

Another  application  of  the  Magic  Wand  Theorem  which is easy to
describe  is  the  following advance in \textit{Illumination Problem}
asking,  whether a room with mirrored walls can always be illuminated
by  a single point light source. R.~Penrose designed in 1958 a planar
room  with  walls made from flat and elliptic mirrors that always has
dark regions no matter where you place a candle in this room. In 1995
G.~Tokarsky  constructed a polygonal room with a similar property: it
has  one  dark point if the idealized candle is placed at the correct
point.  Using  the Magic Wand Theorem, S.~Leli\`evre, T.~Monteil, and
B.~Weiss,   proved   in~\cite{Lelievre:Monteil:Weiss}  that  for  any
translation  surface  $M$,  and any point $x\in M$, the set of points
$y$ which are not illuminated by $x$ is always finite.

One  should  not  have  an  impression  that  the theory developed by
A.~Eskin,  M.~Mirzakhani,  A.~Mohammadi and other researchers in this
area  is designed to serve billiards. A billiard in a polygon is just
a  cute  way  to describe certain class of dynamical systems; we have
seen  that  same  kind  of  dynamical  systems  appear in solid state
physics, in conductivity theory, and so on.

The  result  of  A.~Eskin  and M.~Mirzakhani opens a new way to study
moduli  spaces,  which  in  the last several decades became a central
object both in mathematics and in theoretical physics. We do not know
yet  all  possible applications of the Magic Wand Theorem which might
be obtained in this direction.

The  integral  calculus  was  partly  developed  by Kepler (a century
before  Newton  and  Leibniz)  in order to measure the volume of wine
barrels.  Who  could  imagine  at that time that volume of a solid of
revolution  would  be  discussed  in  any textbook of mathematics for
beginners  and  that  the integral calculus would become an essential
part  of  all  contemporary  engineering.  The theorem proved by Alex
Eskin  and  Maryam  Mirzakhani  is  so  beautiful and powerfull that,
personally,  I have no doubt that it would find numerous applications
far beyond our current imagination.
\smallskip

\textbf{Acknowledgements.}
I    would    like    to    thank   Jon   Chaika,   Alex~Eskin,   and
\mbox{Pascal~Hubert}  for  valuable comments. I am extremely grateful
to Giovanni~Forni, to S\'ebastien Gou\"ezel, to \mbox{Curt~McMullen},
and to my father, Vladimir~Zorich, who carefully read the preliminary
version  of  this  article  and made numerous important remarks which
helped to improve both the content and the presentation.

I thank Maryam~Mirzakhani for the Magic Wand and for the photographs.

\end{document}